\documentclass[12pt]{amsart}
\usepackage[english]{babel}
 \usepackage[utf8x]{inputenc} 
\usepackage{amsmath}
\usepackage{cmll}
\usepackage{amssymb}
\usepackage{amsthm}
\usepackage{xcolor}
\usepackage{mathrsfs}
\usepackage[all]{xy}
\usepackage{indentfirst}
\usepackage{fancyhdr}
\usepackage{newlfont}
\usepackage{setspace}
\usepackage{verbatim}
\usepackage{hyperref} 
\AtBeginDocument{\def\MR#1{}} 

\theoremstyle{definition}

\newtheorem*{theorem*}{Theorem}

\theoremstyle{remark}

\numberwithin{equation}{section}

\newcommand{\qui}{\alpha}
\newcommand{\quo}{\beta}
\newcommand{\qua}{\gamma}

\newcommand{\vertiii}[1]{{\left\vert\kern-0.25ex\left\vert\kern-0.25ex\left\vert #1 
    \right\vert\kern-0.25ex\right\vert\kern-0.25ex\right\vert}}
\newcommand{\norm}[1]{\left\lVert#1\right\rVert}

\begin{document}

\title{A remark on vanishing geodesic distances in infinite dimensions}

\author{Valentino Magnani}
\address{Dipartimento di Matematica, Universit\`a di Pisa}
\curraddr{Dipartimento di Matematica, Universit\`a di Pisa}
\email{valentino.magnani@unipi.it}
\thanks{The first author was supported by the University of Pisa, Project PRA 2018 49.}

\author{Daniele Tiberio}
\curraddr{Dipartimento di Matematica, Universit\`a di Pisa}
\email{danieletiberio@gmail.com}


\date{}



\begin{abstract}
We observe that a vanishing geodesic distance arising from a weak Riemannian metric in a Hilbert manifold can be constructed.
\end{abstract}

\maketitle

\markright{Vanishing geodesic distance in Hilbert manifolds}

It is a well known fact that in a connected and finite dimensional Riemannian manifold taking the infimum among all lengths of curves connecting two points yields a distance.
Understanding whether the analogous procedure still gives a distance for an infinite dimensional manifold is considerably more difficult, when a {\em weak Riemannian metric} is fixed. These metrics are smooth symmetric tensors $g$ on $TM$, with the property that 
\[
g_p(v,v)>0
\]
for every $p\in M$ and $v\in T_pM\setminus\{0\}$. However, it is not required that the associated dual mapping from $T_pM\to T_pM^*$, $v\to  g_p(v,\cdot)$ is an isomorphism of Hilbert spaces. Such additional condition only occurs for {\em strong Riemannian metrics}, see \cite[Definition 5.2.12]{AMR88} for more information. If a manifold $M$ modelled on an infinite dimensional Fr\'echet space $E$ is endowed with a \textit{strong Riemannian metric}, then the model $E$ has a Hilbert space structure and the geodesic distance is actually a distance, see \cite{Klingenberg1995}, \cite{Lang1999Fundamentals} and \cite{Bru2018-Notes} for more information. Clearly strong and weak Riemannian metric coincide on finite dimensional manifolds.  

Given two points in a connected and possibly infinite dimensional manifold $M$ equipped with a weak Riemannian metric, we can clearly define the associated length functional in the usual way. 
For $p,q\in M$ we define the class $\Gamma(p,q)$ of all piecewise smooth curves $\gamma:[0,1]\to M$ such that
$\gamma(0)=p$ and $\gamma(1)=q$. 
If 
\[
L_g(\gamma)=\int_0^1 \sqrt{g_{\gamma(t)}(\dot\gamma(t),\dot\gamma(t))}\,dt
\]
is the standard lengh of a piecewise smooth curve in $M$, we set
\[
d_g(p,q)=\inf\left\{L_g(\gamma):\gamma\in\Gamma(p,q) \right\}.
\]
In infinite dimensional manifolds $d_g$ is a pseudometric. In general it may vanish on distinct points.
Important cases where this phenomenon occurs are diffeomorphism groups and spaces of immersions, that have also interesting applications in shape analysis \cite{BBM14} and computational anatomy \cite{GM98}.
Examples of vanishing geodesic distances have been provided in \cite{EP1993}, \cite{MM05} and \cite{MM06}, see also
\cite{BHP2019preprint} for other recent results. Here the vanishing geodesic distances where constructed in Fr\'echet manifolds. 

It is rather natural to ask whether simple examples of vanishing geodesic distances can be found in Banach or Hilbert manifolds. Our motivations go back to the aim of understanding some aspect of the geometry of homogeneous groups in infinite dimensions. In connection with a Rademacher-type differentiability theorem, some examples of infinite dimensional homogeneous groups have been provided in \cite{MR2014}, using product of spaces of sequences $\ell^p$. These Banach Lie groups can be also equipped with a left invariant distance, that is homogeneous with respect to the groups dilations. We have an additional motivation in understanding whether these Banach Lie groups admit weak Riemannian metrics that give a geodesic distance. More general constructions of infinite dimensional metric Lie groups can be found in \cite{LeDLiMoi2018preprint}.
Many recent works have considered sub-Riemannian manifolds of infinite dimensions under different perspectives. We mention only a few of them, as \cite{GroMarVas2015}, \cite{Arg2016preprint}, \cite{ArgTre2017}, see also references therein. Clearly the list could be enlarged. 

We provide an answer to the above question, by showing a simple example of Hilbert manifold equipped with a weak Riemannian metric whose geodesic distance is everywhere vanishing.

Let $\ell^2$ be the linear space of real-valued and square-summable sequences. We equip $\ell^2$ with the standard scalar product 
$\langle\cdot{}\, , \cdot{}\rangle$, whose norm is ${\norm{x}}=(\sum_{k=1}^\infty{\vert{x_k}\vert}^2)^{1/2}$ for any $x\in\ell^2$. 
Let $A: \ell^2 \to \ell^2$ be the operator that maps $x \in l^2$ to $Ax\in l^2$, defined as
\[
(Ax)_k = \frac{1}{k^4} x_k
\]
for all $k\ge1$. Let $B: \ell^2 \times \ell^2 \to \mathbb{R} $ be the bilinear, symmetric map given 
by $B(x,y) = \langle x,Ay\rangle$.

We consider $\ell^2$ as Hilbert manifold, hence for $p$ in $\ell^2$ and $v,w\in T_p(\ell^2)$, we define the 
weak Riemannian metric 
\[
g_p(v,w)= e^{- \norm{p}^2 } B(v,w),
\]
where have canonically identified the tangent spaces of $\ell^2$ with $\ell^2$ itself.

We will show that for any two distinct points $p$ and $q$ of $\ell^2$ we have $d_g(p,q)=0$. Consider the standard orthonormal basis
${\{e_n\}}_{n=1}^{\infty}$ of $\ell^2$, where $e_1= (1,0,\ldots)$, $e_2 = (0,1,0,\ldots) $ and so on. 
For each positive integer $ n \in \mathbb{N}$, we consider the line segment 
from $p$ to $p + ne_n$, i.e., the curve 
$$
\qui_n(t) = p + tne_n, 
$$ 
where $t\in[0,1]$.   
Then we take the line segment from $p + ne_n$ to $ q + ne_n $ given by
$$\quo_n(t) = p + ne_n + t(q-p)$$
and finally the line segment $\qua_n$ from 
$q+ne_n$ to $q$ given by 
$$\qua_n(t) = q + (1-t)ne_n,$$
where $t$ always varies in $[0,1]$. We join the three curves $\qui_n$, $\quo_n$, and $\qua_n$, then obtaining
a curve $\epsilon_n$ that connects $p$ to $q$. We observe that
$$
\qui'_n(t) = ne_n, 
\,
\quad
\quo'_n(t) = q-p,
\,\quad
\qua'_n(t) =  - ne_n
$$
Our claim follows if we show that $L_g(\qui_n)$, $L_g(\quo_n)$, and  
$L_g(\qua_n)$ converge to zero as $n\to\infty$. Indeed, we get
\begin{align*}
L_g(\qui_n) &= 
\int_0^1 e^{ \frac{- \norm{\qui_n (t)}^2} {2} } \sqrt{B(ne_n,ne_n)}\,dt  \\
&\leq \int_0^1 \sqrt{B(ne_n,ne_n)}\,dt = \sqrt{ \langle ne_n,\frac{1}{n^3}e_n\rangle }=
\frac{1}{n}.
\end{align*}
In the same way, we obtain
$$
L_g(\qua_n) = \int_0^1 e^{ \frac{- \norm{\qua_n (t)}^2} {2} }\sqrt{B(-ne_n,-ne_n)}
\,dt \leq  \int_0^1\sqrt{B(ne_n,ne_n)}
\,dt = 
\frac{1}{n}.
$$
Another simple computation can be carried out for $\quo_n$. We have
\begin{align*}
L_g(\quo_n) &= \int_0^1 e^{ \frac{- \norm{\quo_n (t)}^2} {2} } \sqrt{B(q-p,q-p)} \,dt \\
&\le \sqrt{B(q-p,q-p)} \int_0^1 e^{-\frac{n^2}2+n\|p+t(q-p)\|}\,dt \\
&\leq \sqrt{B(q-p,q-p)} e^{-\frac{n^2}{2} + n (\|p\|+\|q-p\|)}.
\end{align*}
We have shown that $L_g(\epsilon_n)$ converges to zero therefore $d_g(p,q)=0$. 

We have proved the following.

\begin{theorem*}
	There exists a weak Riemannian metric $g$ on $\ell^2$ such that	$d_g\equiv 0$. 
\end{theorem*}

\section*{Acknowledgement} 
The authors wish to thank Erlend Grong for suggesting additional references.

\bibliography{References}
\bibliographystyle{amsalpha}

\end{document}